\input amstex
\documentstyle{amsppt}
\document
\topmatter
\title {On the Cauchy Transform of Weighted Bergman Spaces}
\endtitle
\author {S. Merenkov}
\endauthor
\abstract
The problem of describing the range of a Bergman space $B_2(G)$ under the Cauchy transform $K$ for a Jordan domain $G$ was solved by Napalkov (Jr) and Yulmukhametov [NYu1]. It turned out that $K(B_2(G))=B_2^1(C\overline G)$ if and only if $G$ is a quasidisk; here $B_2^1(C\overline G)$ is a Dirichlet space of the complement of $\overline G$. The description of $K(B_2(G))$ for an integrable Jordan domain is given in [M]. In the present paper we give a description of $K(B_2(G,\omega))$ analogous to the one given in [M] for a weighted Bergman space $B_2(G,\omega)$ with a weight $\omega$ which is constant on level lines of the Green function of $G$. In case $G=\Bbb D$, the unit disk, and under some additional conditions on the weight $\omega$, $K(B_2(\Bbb D,\omega))=B_2^1(C\overline{\Bbb D}, \omega^{-1})$, a weighted analogy of a Dirichlet space. 
\endabstract
\endtopmatter

\head {1. Introduction}
\endhead

Let $G$ be a bounded domain in the complex plane $\Bbb C$ whose boundary is a rectifiable Jordan curve and let $\varphi$ be a conformal map of the unit disk $\Bbb D$ onto $G$. Let $\omega(t)$ be a positive measurable function on $(0,1]$ that is called a weight. We consider weights $\omega$ for which the integral 
$$ 
\iint_G\omega(1-|\psi(z)|)\, dm_2(z) \tag 1.1
$$
converges, where $\psi$ is the inverse function of $\varphi$ and $dm_2$ is the Lebesgue area measure.

Introduce a weighted Bergman space:
$$
B_2(G,\omega)=\left\{ g(z)\in {\text{Hol}}(G), \|g\|_{B_2}=\biggl(\iint_G|g(z)|^2\omega(1-|\psi(z)|)dm_2(z)\biggr)^{\frac12}<\infty\right\}.
$$

Consider the weighted Cauchy transform $K$ of functions from the space $B_2(G,\omega)$:
$$
(Kg)(\zeta)=\frac1{\pi}\iint_G\frac{\overline{g(z)}\omega(1-|\psi(z)|)}{z-\zeta}\, dm_2(z),\qquad \zeta\in C\overline G,\quad g\in B_2(G,\omega), \tag 1.2
$$
where $C\overline G$ means the complement of $\overline G$. The integrability of $\overline{g(z)}\omega(1-|\psi(z)|)$ follows from the convergence of (1.1).

The problem of describing the range of $B_2(G,\omega)$ under the Cauchy transform for $\omega=1$ and different domains was studied in [NYu], [NYu1], [M]. In the present paper we are concerned with studying the same question for weighted spaces. The main result of the paper is a theorem describing the range in terms of spaces $W_{\omega}(0,2\pi)$ which are introduced in $\S$2. The method of proof of the main theorem is similar to the one in [M], but here we use the approximation of functions insted of the domain $G$ when describing $K(B_2(G,\omega))$.

\head {2. More Spaces}
\endhead

To state the theorem of this paper we need to introduce some more spaces, but first we give an example.
\subhead {Example}
\endsubhead
Let $g\in B_2(\Bbb D,\omega),\,\, g(z)=\sum_0^{\infty}a_kz^k$. Taking as a conformal map of $\Bbb D$ onto itself $\varphi(z)= z$ we can easily calculate the norm $\|g\|_{B_2}=(\pi\sum_0^{\infty}|a_k|^2\omega_k)^{1/2}$, where $\{\omega_k=2\int_0^1r^{2k+1}\omega(1-r)\, dr\}_0^{\infty}$ which is a positive decreasing sequence. Also we can compute the Cauchy transform $Kg$:
$$
(Kg)(\zeta)=\frac1{\pi}\iint_{\Bbb D}\frac{\overline{g(z)}\omega(1-|z|)}{z-\zeta}dm_2(z)=\sum_1^{\infty}\frac{b_k}{\zeta^k}, \quad{\text{where}}\,\, b_k=-\overline{a}_{k-1}\omega_{k-1}.
$$
So K is an isometry between $B_2(\Bbb D,\omega)$ and the space
$$
B_2^1(C\overline{\Bbb D},\omega)=\biggl\{\gamma(\zeta)\in {\text{Hol}}(C\overline{\Bbb D}),\,\, \gamma(\zeta)=\sum_1^{\infty}\frac{b_k}{\zeta^k},\,\, \|\gamma\|_{B_2^1}=\biggl(\pi\sum_1^{\infty}\frac{|b_k|^2}{\omega_{k-1}}\biggr)^{\frac12 }<\infty\biggr\}.
$$

Now we introduce the spaces $W_{\omega}(0,2\pi)$ and $\overline{W}_{\omega}(0,2\pi)$:
$$
\aligned
W_{\omega}(0,2\pi)=\biggl\{f(e^{i\theta})\in L^1(0,2\pi),\,\, f(e^{i\theta})\sim\sum_{k=-\infty}^{\infty}f_ke^{ik\theta},\\
\rho(f)=\biggl(\pi\sum_{k=1}^{\infty}\frac{|f_{-k}|^2}{\omega_{k-1}}\biggr)^{\frac12}<\infty\biggr\}.   
\endaligned
\tag 2.1
$$

Let $\{G_n\}_1^{\infty}$ be a sequence of domains approximating the domain $G$ from the inside, that is a) $\overline{G}_n\subset G_{n+1}\subset G,\,\, n\in \Bbb N$; b) $\cup_{n=1}^{\infty}G_n=G$. We say that a function $\gamma(\zeta)$ analytic in $C\overline G$ belongs to $\overline W_{\omega}(0,2\pi)$ if for some sequence of domains $\{G_n\}_1^{\infty}$ approximating the domain $G$ from the inside there exists a sequence of functions $\{\gamma_n\}_1^{\infty}$ satisfying the following conditions:
\roster
\item $\gamma_n$ is analytic in $C\overline G_n$;
\item $\gamma_n$ converges to $\gamma$ uniformly on every compact $K\subset C\overline G$;
\item  $\sup_n\rho(\gamma_n\circ\varphi)<\infty$, where $\circ$ means composition.
\endroster

The space $\overline W_{\omega}(0,2\pi)$ is a kind of closure of the system of functions $\gamma$ analytic in $CG$ with $\gamma\circ\varphi\in W_{\omega}(0,2\pi)$, but the closure is in the topology of uniform convergence on compact sets of $C\overline G$.

In order to prove the theorem of this paper we need to impose an additional condition on the weight $\omega$. Let $P_2(G,\omega)$ be the closure of all analytic polynomials in  $B_2(G,\omega)$. Throughout the paper we assume that $\omega$ is chosen in such a way that 
$$
P_2(G,\omega)=B_2(G,\omega).  \tag 2.2
$$
In connection with (2.2) see, for example, [V], [B], [B1]. One condition that guarantees that $\omega$ satisfies (2.2) is (see [V]) 
$$
\aligned
\cases
t\log1/{\omega(t)}\uparrow+\infty,\quad {\text{as}} t\downarrow 0,\\
\int_0\log\log1/{\omega(t)}dt=\infty.
\endcases
\endaligned
$$

\head {3. The Main Theorem}
\endhead

\proclaim{Theorem}
Let $G$ be a bounded domain in $\Bbb C$ whose boundary is a rectifiable Jordan curve, and let $\omega(t)$ be a positive measurable function on $(0,1]$ such that integral (1.1) converges and condition (2.2) is satisfied. Then a function $\gamma$ analytic in $C\overline G$ belongs to $K(B_2(G,\omega))$ if and only if $\gamma \in\overline W_{\omega}(0,2\pi)$, i.e.
$$
K(B_2(G,\omega))=\overline W_{\omega}(0,2\pi).
$$ 
\endproclaim
\demo{Proof}
Suppose that $\gamma\in\overline W_{\omega}(0,2\pi)$. This means that there exists a sequence of domains $\{G_n\}_1^{\infty}$ approximating the domain $G$ from the inside and a sequence of functions $\{\gamma_n\}_1^{\infty}$ with
\roster
\item$\gamma_n$ analytic in $C\overline G_n$;
\item$\gamma_n\rightarrow\gamma$ uniformly on every compact set $K\subset C\overline G$;
\item$\sup_n\rho(\gamma_n\circ\varphi)<\infty$, where $\rho$ is from (2.1).
\endroster

For every $n\in\Bbb N$ consider the linear functional $\Bbb F_n$ on functions $h\in{\text{Hol}}(\overline G)$ given by 
$$\Bbb F_n(h)=-\frac1{2\pi i}\int_{\partial\, G}\gamma_n(\xi)h(\xi)\, d\xi. \tag 3.1
$$
We are going to prove that $\Bbb F_n$ is a bounded linear functional in the norm of the space $B_2(G,\omega)$.
If we make the change of variables $\xi=\varphi(e^{i\theta})$ in the integral in (3.1) we get
$$
-\frac1{2\pi i}\int_0^{2\pi}\gamma_n(\varphi(e^{i\theta}))h(\varphi(e^{i\theta}))\varphi_{\theta}'(e^{i\theta})\, d\theta.
$$
Since the boundary $\partial\, G$ is a rectifiable curve, the function $h(\varphi(e^{i\theta}))\varphi_{\theta}'(e^{i\theta})$ is the restriction of the function $h_1(z)zi$ on the unit circumference, where $h_1(z)=h(\varphi(z))\varphi'(z)$ [G, p. 405]. Thus, the last integral is equal to
$$
-\frac1{2\pi}\int_0^{2\pi}\gamma_n(\varphi(e^{i\theta}))h_1(e^{i\theta})e^{i\theta}d\theta.
$$
The operator $\Phi:\,\, h\rightarrow h_1$ is an isometry between the spaces $B_2(G,\omega)$ and $B_2(\Bbb D,\omega)$, that is 
$$
\|h_1\|_{B_2}=\|h\|_{B_2}. \tag 3.2
$$
Let $\gamma_n(\varphi(e^{i\theta}))\sim \sum_{k=-\infty}^{\infty}c_k^ne^{ik\theta},\,\, h_1(z)=\sum_{k=0}^{\infty}a_kz^k$, then
$$
-\frac1{2\pi i}\int_{\partial\, G}\gamma_n(\xi)h(\xi)\, d\xi=-\sum_{k=0}^{\infty}c_{-(k+1)}^na_k.
$$
Applying the Cauchy-Schwarz inequality we get
$$
|\Bbb F_n(h)|
\leq \biggl(\sum_{k=1}^{\infty}\frac{|c_{-k}^n|^2}{\omega_{k-1}}\biggr)^{\frac12}\biggl(\sum_{k=0}^{\infty}\omega_k|a_k|^2\biggr)^{\frac12}=\frac1{\pi}\rho(\gamma_n\circ\varphi)\|h\|_{B_2}, \tag 3.3
$$
where we used (3.2). 
Since polynomials are dense in $B_2(G,\omega)$ the functional $\Bbb F_n$ is uniquely extended to a bounded linear functional, which we also call $\Bbb F_n$, on $B_2(G,\omega)$ with $\|\Bbb F_n\|\leq\frac1{\pi}\rho(\gamma_n\circ\varphi)$. 

The space $B_2(G,\omega)$ is a Hilbert space, hence for every $n\in\Bbb N$ there exists a function $g_n\in B_2(G,\omega)$, such that
$$
\Bbb F_n(h)=\frac1{\pi}\iint_{G}h(z)\overline{g_n(z)}\omega(1-|\psi(z)|)\, dm_2(z).
$$
Moreover, 
$$
\|g_n\|_{B_2}=\|\Bbb F_n\|. \tag 3.4
$$
From (3.3), (3.4), condition (3) for the sequence $\{\gamma_n\}_1^{\infty}$ and applying the Banach-Alaoglu theorem we conclude that there exists a bounded linear functional $\Bbb F$ on $B_2(G,\omega)$ such that some subsequence $\{\Bbb F_{n(k)}\}_{k=1}^{\infty}$ of $\{\Bbb F_n\}_{n=1}^{\infty}$ converges to $\Bbb F$ in the weak-star topology. Let $g\in B_2(G,\omega)$ be a function such that
$$
\Bbb F(h)=\frac1{\pi}\iint_{G}h(z)\overline{g(z)}\omega(1-|\psi(z)|)dm_2(z),\quad h\in B_2(G,\omega)
$$
and $\|g\|_{B_2}=\|\Bbb F\|$. Computing the value of $\Bbb F_n$ at $h(z)=1/(z-\zeta), \, \zeta\in C\overline G$ we get
$$
\gamma_n(\zeta)=\Bbb F_n\biggl(\frac1{z-\zeta}\biggr)=\frac1{\pi}\iint_G\frac{\overline{g_n(z)}}{z-\zeta}\omega(1-|\psi(z)|)dm_2(z).
$$
Letting $n=n(k)$ tend to infinity and using condition (2) for the sequence $\{\gamma_n\}_1^{\infty}$ reveals that 
$$
\gamma(\zeta)=\frac1{\pi}\iint_G\frac{\overline{g(z)}\omega(1-|\psi(z)|)}{z-\zeta}dm_2(z), \quad \zeta\in C\overline G,
$$
i.e. $\gamma(\zeta)=(Kg)(\zeta)$.

To prove the converse, take any function $\gamma\in K(B_2(G,\omega))$, i.e., by (1.2),
$$
\gamma(\zeta)=\frac1{\pi}\iint\frac{\overline{g(z)}\omega(1-|\psi(z)|)}{z-\zeta}dm_2(z),\quad \zeta\in C\overline G,
$$
$ g\in B_2(G,\omega)$. We shall now prove that $\gamma$ belongs to $\overline W_{\omega}(0,2\pi)$.

For any $n\in\Bbb N$ consider a function $\alpha_n$ such that
\roster
\item $\alpha_n$ is continuous on $[0,1],\,\, 0\leq\alpha_n(t)\leq1$;
\item $\alpha_n(t)=0,\,\, t\in[0,1/n];\,\, \alpha_n(t)=1,\,\, t\in[2/n,1]$.
\endroster
Set 
$$
\gamma_n(\zeta)=\frac1{\pi}\iint_G\frac{\overline{g(z)}\omega(1-|\psi(z)|)}{z-\zeta}\alpha_n(1-|\psi(z)|)\, dm_2(z). \tag 3.5
$$
First, it is evident that if we take $G_n=\{z\in G:\,\, |\psi(z)|<1-1/n\}$ then $\gamma_n$ is analytic in $C\overline G_n$. Next we prove that $\gamma_n\rightarrow \gamma$ uniformly on a compact $K\subset C\overline G$.
$$
\aligned
|\gamma(\zeta)-\gamma_n(\zeta)|&=\frac1{\pi}\biggl|\iint_G\frac{\overline{g(z)}\omega(1-|\psi(z)|)}{z-\zeta}(1-\alpha_n(1-|\psi(z)|))\, dm_2(z)\biggr|\\
&\leq\frac{C_K}{\pi}\|g\|_{B_2}\biggl(\iint_G(1-\alpha_n(1-|\psi(z)|))^2\omega(1-|\psi(z)|)\, dm_2(z)\biggr)^{\frac12},
\endaligned
$$
where the constant $C_K$ depends on the compact set $K$.
The last integral, after changing variables, becomes 
$$
\iint_{1-2/n\leq|z|<1}|\varphi'(z)|^2\omega(1-|z|)\, dm_2(z).
$$
Hence $\sup_K|\gamma(\zeta)-\gamma_n(\zeta)|\rightarrow0$ as $n\rightarrow\infty$.

Finally, it remains to prove that $\sup_n\rho(\gamma_n\circ\varphi)<+\infty$, where $\gamma_n$ is defined by (3.5). For $f\in L^1(\partial\,\Bbb D),\,\, f(e^{i\theta})\sim\sum_{-\infty}^{\infty}f_ke^{ik\theta}$ the corresponding Cauchy type integral
$$
F(\zeta)=\frac1{2\pi i}\int_{\partial\, \Bbb D}\frac{f(t)}{t-\zeta}dt,\quad \zeta\in C\overline{\Bbb D}
$$
has the Taylor expansion $F(\zeta)=\sum_{k=1}^{\infty}f_{-k}/\zeta^k,\,\, \zeta\in C\overline{\Bbb D}$. With this in mind, consider
$$
F_n(\zeta)=\frac1{2\pi i}\int_{\partial\,\Bbb D}\frac{\gamma_n\circ\varphi(t)}{t-\zeta}dt,\quad \zeta\in C\overline{\Bbb D}. \tag 3.6
$$
If we substitute (3.5) in (3.6), using $\zeta=\varphi(t),\,\, t\in \Bbb D$, change the order of integration, and compute the inner integral, we get
$$
F_n(\zeta)=-\frac1{\pi}\iint_G\frac{\overline{g(z)}\omega(1-|\psi(z)|)}{(\psi(z)-\zeta)\varphi'(\psi(z))}\alpha_n(1-|\psi(z)|)\, dm_2(z).
$$
Make the change of variable $w=\psi(z)$ to obtain
$$
\aligned
F_n(\zeta)&=-\frac1{\pi}\iint_{\Bbb D}\frac{\overline{g(\varphi(w))\varphi'(w)}\omega(1-|w|)}{w-\zeta}\alpha_n(1-|w|)\, dm_2(w)\\
&=\frac1{\pi}\sum_{k=0}^{\infty}\frac1{\zeta^{k+1}}\iint_{\Bbb D}\overline{g(\varphi(w))\varphi'(w)}w^k\omega(1-|w|)\alpha_n(1-|w|)\, dm_2(w).
\endaligned
$$

As was noted above, $g(\varphi(z))\varphi'(z)\in B_2(\Bbb D,\omega)$. Let $g(\varphi(z))\varphi'(z)=\sum_{k=0}^{\infty}a_jz^j$. Then $F_n(\zeta)=\sum_{k=1}^{\infty}b_k^n/\zeta^k$, where $b_k^n=\overline a_{k-1}2\int_0^1r^{2k-1}(\alpha_n\omega)(1-r)\, dr$. We need to check that $\sup_n(\pi\sum_{k=!}^{\infty}|b_k^n|^2/\omega_{k-1})^\frac12<\infty$. Since $0\leq\alpha_n\leq1$,
$$
\frac{|b_k^n|^2}{\omega_{k-1}}\leq |a_{k-1}|^2\omega_{k-1}.
$$
Hence
$$
\sup_n\biggl(\pi\sum_{k=1}^{\infty}\frac{|b_k^n|^2}{\omega_{k-1}}\biggr)^{\frac12}\leq (\pi\sum_{k=0}^{\infty}|a_k|^2\omega_k)^{\frac12}=\|g\|_{B_2}.
$$

Thus the theorem is proved.
 
\enddemo

\head{4. Description  of $W_{\omega}(0,2\pi)$.}
\endhead

Let $f\in L^1(\partial\,\Bbb D)$, that is $f(e^{i\theta})\in L^1(0,2\pi),\,\, f(e^{i\theta})\sim \sum_{k=-\infty}^{\infty}f_ke^{ik\theta}$. The Cauchy-type integral $F$ corresponding to $f$ is $F(\zeta)=-\sum_{k=1}^{\infty}f_{-k}/\zeta^k,\,\, \zeta\in C\overline{\Bbb D}$. We prove that under some condition on the weight $\omega,\,\, f\in W_{\omega}(0,2\pi)$ if and only if $F\in B_2^1(C\overline{\Bbb D},\omega^{-1})$, where
$$
\aligned
B_2^1(C\overline{\Bbb D},\omega^{-1})=\biggl\{f(\zeta)\in{\text {Hol}}(C\overline{\Bbb D}), F(\infty)&=0,\\ 
\|F\|_{B_2^1}&=\biggl(\iint_{C\overline{\Bbb D}}|F'(\zeta)|^2\frac{dm_2(\zeta)}{\omega(1-1/|\zeta|)}\biggr)^{\frac12}\biggr\}.
\endaligned
$$

\proclaim{Proposition}
Assume that $\omega$ satisfies the following condition:
$$\sup_k\biggl\{(k+1)^2\int_0^1r^{2k+1}\omega(1-r)\, dr\int_0^1r^{2k+1}\frac{dr}{\omega(1-r)}\biggr\}\leq C<\infty. \tag 4.1
$$
Then $f \in W_{\omega}(0,2\pi)$ if and only if $F\in B_2^1(C\overline{\Bbb D},\omega^{-1})$.
\endproclaim
\demo{Proof}

First we show that 
$$
\inf_k\biggl\{(k+1)^2\int_0^1r^{2k+1}\omega(1-r)\, dr\int_0^1r^{2k+1}\frac{dr}{\omega(1-r)}\biggr\}\geq c, \tag 4.2
$$
where $c$ is some positive constant:
$$
\frac1{4(k+1)^2}=\biggl(\int_0^1r^{2k+1}dr\biggr)^2\leq\int_0^1r^{2k+1}\omega(1-r)\, dr\int_0^1r^{2k+1}\frac{dr}{\omega(1-r)},
$$
which is (4.2) with $c=1/4$.
Now we are ready to prove the proposition.
$$
\aligned
\iint_{C\overline{\Bbb D}}|F'(\zeta)|^2\frac{dm_2(\zeta)}{\omega(1-1/|\zeta|)}=2\pi\sum_{k=1}^{\infty}k^2|f_{-k}|^2\int_1^{\infty}\frac{dr}{r^{2k+1}\omega(1-1/r)}\\
=2\pi\sum_{k=1}^{\infty}k^2|f_{-k}|^2\int_0^1\frac{r^{2k-1}dr}{\omega(1-r)}\asymp\pi\sum_{k=1}^{\infty}\frac{|f_{-k}|^2}{\omega_{k-1}},
\endaligned
$$
where $a\asymp b$ means that there exist two positive constants $m, M$ such that $ma\leq b\leq Ma$; here we used (4.1), (4.2).
This proves the proposition. 
\enddemo
\proclaim{Remark}
The condition (4.1) is similar to the Mackenhaupt condition [Ga, p. 254].
\endproclaim

\enddocument